\documentclass[10pt]{article}
\usepackage{amsfonts, epsfig, amsmath, amssymb, color,mathabx}

\usepackage[english]{babel}

\renewcommand {\epsilon}{\varepsilon}

\newtheorem{thm}{Theorem}[section]

\newtheorem{lem}{Lemma}[section]

\newtheorem{exm}{Example}[section]

\DeclareMathSymbol{\ophi}{\mathalpha}{letters}{"1E}

\renewcommand{\phi}{\varphi}

\newcommand{\be}{\begin{equation}}
\newcommand{\ee}{\end{equation}}
\newcommand{\ben}{\begin{equation*}}
\newcommand{\een}{\end{equation*}}

\newcommand{\ba}{\begin{equation}\begin{aligned}}
\newcommand{\ea}{\end{aligned}\end{equation}}

\newenvironment{proof}{\par\noindent{\bf Proof:}}{\hfill$\blacksquare$\par}

\newfont{\cyrfnt}{wncyr10}
\def\J3{\cyrfnt{\rm \u{\cyrfnt I}}}
\def\j3{\cyrfnt{\rm \u{\cyrfnt i}}}

\allowdisplaybreaks[4]

\begin{document}

\title{FIRST EXIT TIMES OF SOLUTIONS OF STOCHASTIC DIFFERENTIAL EQUATIONS 
DRIVEN BY MULTIPLICATIVE L\'EVY NOISE WITH HEAVY TAILS\footnote{Stochastics and Dynamics, 
Vol.\ 11, Nos.\ 2 \& 3 (2011) 495--519} 
}

\author{ILYA PAVLYUKEVICH\\ \\
Institut f\"ur Stochastik\\ Fakult\"at f\"ur Mathematik und Informatik\\
Friedrich--Schiller--Universit\"at Jena\\ Ernst--Abbe--Platz 2, Jena\\
07743, Germany\\
\texttt{ilya.pavlyukevich@uni-jena.de}}


%

\date{}

\maketitle
%

\begin{center}
\begin{small}\textit{Dedicated to Peter Imkeller on the occasion of his 60th birthday, \\with friendship and respect\\\null}\end{small}
\end{center}
\begin{abstract}
In this paper we study first exit times from a bounded domain of a gradient dynamical system 
$\dot Y_t=-\nabla U(Y_t)$ perturbed by a small multiplicative 
L\'evy noise with heavy tails.  A special attention is paid to the 
way the multiplicative noise is introduced. In particular we determine the asymptotics of the first exit time
of solutions of It\^o, Stratonovich and Marcus canonical SDEs.
\end{abstract}



\numberwithin{equation}{section}
\section{Introduction}

In many models of natural phenomena the state of a system is 
described by a deterministic ordinary differential equation of the form
\begin{equation}
\label{eq:y}
Y_t=y+\int_0^t B(Y_s)\, ds,\quad t\geq 0.
\end{equation}
It is often supposed that the vector field $B$ is determined by a function $U$, so that $B=-\nabla U$.
The function $U$ could  be called a climatic pseudo-potential 
\cite{BenziPSV-83,Ditlevsen-99a,Nicolis-82} in geosciences, energy potential in physics 
\cite{Gardiner-04,Kramers-40}
or profit or cost function in economics and optimization \cite{Tu-94}.  
The potential $U$ is often supposed to have several local minima corresponding to the steady states of the system $Y$.
The state space can be decomposed into a number of domains of attraction, so that a solution $Y_t(y)$ 
cannot pass from one domain to the other.

In order to make the models more realistic and allow transitions between the stable states, the system \eqref{eq:y}
is being perturbed by a small random noise, so that  \eqref{eq:y} turns to a random  equation with a small 
parameter.
 Clearly, the properties of the new random system depend on the interplay between the type of the noisy perturbation and the 
underlying deterministic vector field $B$.

Noise can be included into the system in different ways. If the perturbation does not depend on the state
of the system, one usually speaks about additive noise. If the amplitude of the noise depends on the 
state of the system, one speaks about multiplicative perturbations. 

Let for example $Z$ be a regular random process, 
say, with smooth paths, and $F$ be a smooth bounded function. Then 
the perturbed system with multiplicative noise is described by the random ordinary integral equation
\begin{equation}
\label{eq:w}
X_t=x-\int_0^t \nabla U (X_s)\, ds+\varepsilon  \int_0^t F(X_{s})\, dZ_s,\quad t\geq 0,
\end{equation}
where the last integral is understood in Lebesgue--Stieltjes sense and a positive small parameter $\varepsilon$ determines 
the noise amplitude. The situation becomes more complicated if one considers irregular 
perturbations, for instance when $Z$ is a Brownian motion. In this case, the differential $dZ$
is usually understood in the sense of the stochastic It\^o calculus.

There is a lot of literature devoted to the small noise equation \eqref{eq:w}, both from the 
point of view of Mathematics and applications.
The main reference on the large deviations theory and asymptotics of the exit times of equation \eqref{eq:w} 
driven by the Brownian motion $Z$ is 
Freidlin and Wentzell \cite{FreidlinW-98}. In this case, the first exit time of $X$ from a domain around the steady state of 
the underlying deterministic system appears to be exponentially large of the order
$e^{C/\varepsilon^2}$ with the rate $C>0$ being interpreted as the energy the Brownian particle should have in order to
reach the boundary of the domain of attraction. A good exposition of small noise properties of Gaussian SDEs 
with applications can be found in Olivieri and Vares \cite{OlivieriV-03} and Schuss \cite{Schuss-10}. 
Very exact asymptotics of the mean first exit time in the Gaussian case was obtained in Bovier \textit{et al.} 
\cite{BovierEGK-04,BovierGK-05}.

Recently dynamical systems perturbed by small \textit{jump noise with heavy tails} attracted the attention of 
the physical and mathematical community. The physicists' research focuses on the models 
incorporating $\alpha$-stable non-Gaussian L\'evy processes, often referred to as \textit{L\'evy flights}.
Thus Ditlevsen \cite{Ditlevsen-99b,Ditlevsen-99a} proposed an
interesting conjecture about the $\alpha$-stable noise signal in the Greenland ice-core data (see also Hein \textit{et al.}
\cite{HeinImkPavl-09}
on the statistical treatment of this time series). An enhanced, certainly non-exhaustive list of physical references on 
the first exit problem of L\'evy-driven SDEs with stable noises includes Chechkin \textit{et al.} 
\cite{ChechkinSMK-07,ChechkinGKM-05} and Dybiec \textit{et al.} \cite{DybGudHng06,DybGudHng07}.

The mathematical theory of large deviations for general Markov processes can be found in Wentzell \cite{Wentzell-90}.   
To our knowledge, in \cite{Wentzell-90} and in Godovanchuk \cite{Godovanchuk-82} the asymptotic behaviour of the dynamical 
systems with heavy power tails was considered for the first time. 
Opposite to the Gaussian case, the behaviour of such systems is mainly governed by big jumps. Thus, the exit from the 
domain occurs with the help of an only big jump, and the mean exit time does not depend on  
the energy landscape of the underlying dynamical system, but rather on the geometric layout of the stable states and domains of attraction. 
Fine small noise asymptotics of the SDEs with additive heavy tail L\'evy noise 
and their metastable behaviour was studied in Imkeller and Pavlyukevich \cite{ImkellerP-08}. The case of light,
sub- and super-exponential jumps was considered in  Imkeller \textit{et al.} \cite{ImkPavlWetz-09}.
In his very recent work, H\"ogele \cite{Hoegele-Diss} studied the first exit problem and metastability properties 
of solutions 
of the infinite-dimensional stochastic Chafee--Infante equation driven by small heavy tail L\'evy noise.
 
Coming back to the equation \eqref{eq:w} with multiplicative noise, it is necessary to note that the stochastic integral
w.r.t.\ $Z$ allows interpretations different from the It\^o definition, in particular one can consider 
Stratonovich integrals often denoted by $\circ \,dZ$. Even for continuous integrators $Z$,
an interesting question arises, namely, which integral fits a specific real world phenomenon, see 
Arnold \cite{Arnold-74}, Turelli \cite{Turelli77}, van Kampen \cite{vanKampen81}, 
Sethi and Lehoczky \cite{SetLeh81}, Smythe \textit{et al.} \cite{SmyMosMcCCla83} and Sokolov \cite{Sokolov-10} for discussion.
Roughly speaking, It\^o SDEs appear naturally as a continuous approximation of a discrete system, for instance
in financial mathematics or biology. Due to their nice mathematical properties, they are also 
the most popular tools in analysis. On the other hand, Stratonovich SDEs w.r.t.\ continuous integrators $Z$ 
arise naturally as a mathematical idealization of dynamical systems perturbed by regular stochastic processes, 
which takes place in engineering and physical sciences.
Moreover, Stratonovich integrals enjoy a conventional Newton--Leibniz change of variables formula; 
they are also indispensable for constructing SDEs on manifolds.

If the integrator $Z$ is a jump process, for instance an $\alpha$-stable L\'evy process, the simple Newton--Leibniz 
change of variables formula does not hold any longer even for
the Stratonovich integral. To correct this situation, the so called canonical SDEs were introduced be S.\ I.\ Marcus
in \cite{Marcus-78,Marcus-81}. 

In this paper we study multidimensional SDEs of the type \eqref{eq:w}. The random process $Z$ is supposed to be 
a multivariate L\'evy noise with regularly varying (heavy) tails, and the stochastic differential equation will be understood 
in the senses of It\^o, Stratonovich and Marcus. Our study treats the first exit time of the 
perturbed system from a bounded domain around a stable attractor
of the underlying deterministic dynamical system $Y$ in the limit of small noise.

\section{Object of study}

\subsection{The underlying dynamical system} 
We start with a $n$-dimensional gradient system generated by a vector field $-\nabla U$,
\begin{equation*}
Y_t=y-\int_0^t \nabla U(Y_s)\, ds,\quad t\geq 0.
\end{equation*}
We assume that the potential $U$ is a $C^2(\mathbb{R}^n,\mathbb{R})$ function
with globally Lipschitz continuous first derivatives $\partial_i U(x)$, $1\leq i\leq n$, and 
bounded second derivatives $\partial_i\partial_j U(x)$, $1\leq i,j\leq n$.
We also assume that the potential $U$ has a unique global minimum at the origin, $U(0)=0$, 
that is $\nabla U(0)=0$ and the Hesse matrix $(\partial_i\partial_j U(0))_{i,j=1}^n$ is positive definite.

Let $\mathcal{G}\subset \mathbb{R}^n$ be a bounded domain with piece-wise smooth boundary $\partial \mathcal{G}$ such that $0\in \mathcal{G}$. Assume that 
$\langle n(y), -\nabla U(y)\rangle\leq -\delta$ for $y\in \partial\mathcal{G}$ and some $\delta>0$, where $n(y)$ is a unit 
outward normal at $y\in \partial \mathcal{G}$. 

Under these assumptions $0$ is the unique asymptotically stable attractor of the dynamical system $Y_t(y)$, 
$Y_t(y)\to 0$, $t\to\infty$;  for all $y\in \mathcal{G}$ the trajectories $Y_t(y)$ do not 
leave the domain $\mathcal{G}$. 

Finally let 
$F(x)=(F_{ij}(x))_{i,j=1}^{n,m}$, $x\in\mathbb{R}^n$, be a $n\times m$ matrix of smooth bounded real functions with 
Lipschitz continuous bounded  derivatives. Let  $\varepsilon>0$ be a small parameter.

\subsection{ The driving L\'evy process}
On a filtered probability space $(\Omega,\mathcal{F},\mathbf{P})$ satisfying the usual hypothesis we consider an $m$-dimensional L\'evy
process $Z=(Z^1,\dots, Z^m)$ with the characteristic triplet
$(A, \nu,\mu)$ with a non-negative definite $m\times m$ matrix $A$, a vector $\mu\in\mathbb{R}^m$, and a 
L\'evy measure $\nu$ with $\nu(\{0\})=0$ and $\int (1\wedge \|y\|^2)\, \nu(dy)<\infty$.  
In other words 
the characteristic function of $Z$ is given by the L\'evy--Khintchine formula
\begin{equation*}
\mathbf{E} e^{i\langle \lambda, Z_t\rangle}=\exp\Big(
it\langle\lambda,\mu \rangle - t\frac{\langle A\lambda,\lambda\rangle}{2}+t\int (e^{i\langle\lambda, y\rangle}-1-
i\langle\lambda, y\rangle
\mathbb{I}_{\{\|z\|\leq 1\}})\, \nu(dy)
\Big), \lambda\in\mathbb{R}^m.
\end{equation*}
There is also a canonical L\'evy--It\^o representation of $Z$ as a sum
\begin{equation*}
Z_t=\mu t+ W_t+\int_{(0, t]}\int_{0<\|z\|<1} z (N(ds, dz)-ds\nu(dy))+\int_{(0,t]}\int_{\|z\|\geq 1} z N(ds, dz),
\end{equation*}
with $W$ being a Brownian motion with the covariance matrix $A$ and $N$ being a Poisson random measure with the 
intensity measure $\nu$.

To specify the heavy tail property of $Z$ we assume that $\nu$ 
is a regularly varying jump measure at $\infty$. 
Let $H(u)$ denote its tail,
\begin{equation*}
H(u):=\nu(\{z\in\mathbb{R}^m:\ \|z\|\geq u\}).
\end{equation*} 
Then for any $a>0$ the measure $\nu$ enjoys the following 
scaling property:
there is a non-zero Radon measure $m$ on 
$\mathcal{B}(\overline{\mathbb{R}}^m\backslash\{0\})$ with $m(\overline{\mathbb{R}}^m\backslash \mathbb{R}^m)=0$ so that
for any Borel set $A$ bounded away from the origin, $0\notin \overline{A}$, with $m(\partial A)=0$ the relation
\begin{equation*}
m(aA)=\lim_{u\to +\infty}\frac{\nu(auA)}{H(u)}=\frac{1}{a^r}\lim_{u\to +\infty}\frac{\nu(uA)}{H(u)}
=\frac{1}{a^r}m(A)
\end{equation*}
holds for some $r>0$. 
In particular, $H(u)$ is regularly varying at infinity with index $-r$, that is $H(u)=u^{-r}l(u)$ for some positive
slowly varying function $l$. 
The homogeneity property of the limit measure $m$ implies that $m$
assigns no mass to spheres centred at the origin on $\mathbb{R}^m$ and has no
atoms.

For more information on multivariate heavy tails and regular variation we refer the reader to 
Resnick \cite{Resnick-04} and Hult and Lindskog \cite{HultL-06,HultL-06-1}.

\subsection{SDE with multiplicative noise}

In this section we briefly remind the main properties of the  It\^o, Stratonovich and Marcus (canonical) SDEs.

\subsubsection{It\^o SDE}

For simplicity we start in the one-dimensional setting.
Let  $g_t$ be a c\`adl\`ag adapted stochastic process. Then its left-continuous modification $g_{t-}$ 
 is predictable and can be approximated w.r.t.\ a u.c.p.\ topology by simple 
predictable processes $g^{(k)}$ of the form
\begin{equation*}
g^{(k)}_t=g_0\mathbb{I}_{\{0\}}(t)+\sum_{j=1}^k g_{j}\mathbb{I}_{(\tau_j,\tau_{j+1}]}(t),
\end{equation*}
where $0=\tau_0<\cdots<\tau_k$ are stopping times and $g_j^{(k)} $ are $\mathcal{F}_{\tau_j}$ measurable and bounded.
For a L\'evy process $Z$ (or even a semimartingale),
the It\^o stochastic integral of $g$ w.r.t.\ $Z$ is then defined as a limit 
\begin{equation*}
\int_0^t g_{s-}\, dZ_s
:=\lim_{k\to\infty}\sum_{j=0}^k g_j(Z_{\tau_j^{(k)}\wedge t}
-Z_{\tau_{j-1}^{(k)}\wedge t})
\end{equation*}
in the sense of the u.c.p.\ topology, see Chapter II in Protter \cite{Protter-04}.

In particular, one can approximate the It\^o integral by 
\textit{non-anticipating} Riemannian sums.
Indeed, consider a sequence of random partitions
$\tau^{(n)}=\{0=\tau_0^{(n)}\leq \tau_1^{(n)}\leq \cdots\leq \tau_{k_n}^{(n)}<\infty\}$ with 
$\limsup_n \tau_{k_n}^{(n)}=\infty$ a.s., and $\|\tau^{(n)}  \|:=\sup_k| \tau_k^{(n)} -\tau_{k-1}^{(n)}|\to 0$ a.s. Then
\begin{equation*}
\int_0^t g_{s-}\, dZ_s=\lim_{k\to\infty}\sum_{j=0}^k g_{\tau_{j-1}}(Z_{\tau_j^{(k)}\wedge t}
-Z_{\tau_{j-1}^{(k)}\wedge t})
\end{equation*}
in the sense of the u.c.p.\ topology, see Theorem II.21 in Protter \cite{Protter-04}. 
We refer the reader to Applebaum
\cite{Applebaum-09}
and Kunita \cite{Kunita-04} for the theory of stochastic integration w.r.t.\ L\'evy processes, and also to 
Protter \cite{Protter-04} for the general semimartingale theory. 

Now we introduce the It\^o stochastic differential equation with small multiplicative noise. The matrix valued function $F$ given,
we perturb the equation \eqref{eq:y} to obtain
\begin{equation}
\label{eq:i}
X_t=x-\int_0^t \nabla U(X_s)\, ds+\varepsilon \int_0^t F(X_{s-})\, dZ_s.
\end{equation}
In the coordinate form this equation reads
\begin{equation*}
X^i_t=x_i-\int_0^t \partial_{i}U (X_s)\, ds+\varepsilon\sum_{j=1}^m \int_0^t F_{ij}(X_{s-})\, dZ^j_s, \quad 1\leq i\leq n.
\end{equation*}
In particular, under above conditions, there exists a strong solution to the equation \eqref{eq:i},
which is a c\`adl\`ag semimartingale and a strong Markov process,  \cite{Applebaum-09,Kunita-04,Protter-04}.
 
Let $f\in C^2(\mathbb{R}^n,\mathbb{R})$ and $X$ be a the solution of 
\eqref{eq:i}. Then the following change of variables formula (It\^o's formula) holds (Theorem II.33 in Protter \cite{Protter-04}):
\begin{equation*}
\begin{aligned}
f(X_t)=f(x)&+\sum_{i=1}^n \int_0^t \partial_{i} f(X_{s-})\, dX^i_s
+\frac{1}{2}\sum_{i,j=1}^n\int_0^t \partial _i\partial_j f(X_s)\,d[X^i, X^j]^c_s\\
&+
\sum_{s\leq t}\Big(f(X_s)-f(X_{s-})-\sum_{j=1}^n \partial_{i} f(X_{s-})\Delta X^i_s  \Big)
\end{aligned}
\end{equation*}
with $[X^i, X^j]^c$ being the path-by-path continuous part of the quadratic covariation of $X^i$ and $X^j$.

\subsubsection{Stratonovich SDE }

Let again $g_t$ be a c\`adl\`ag adapted stochastic process and  $Z$ be a L\'evy process, such that  the quadratic 
covariation
$[g, Z]$ exists. The Stratonovich integral of $g_{t -}$ w.r.t.\ $Z$ is defined with the help of the It\^o integral as
\begin{equation*}
\int_0^t g_{s-}\circ dZ_s
= \int_0^t g_{s-} \, dZ_s
+\frac{1}{2}[g,Z]^c_t.
\end{equation*} 
The Stratonovich integral can be also interpreted as a limit of Riemannian sums. Let $g$ and $Z$ have no jumps
in common, that is $\sum_{s\leq t}\Delta g_s\Delta Z_s=0$ for all $t\geq 0$, then
 for any sequence of random partitions
$\tau^{(n)}$ we have
\begin{equation*}
\lim_{n\to\infty} \sum_{j=1}^n\frac{g_{\tau_j^{(n)}}+g_{\tau_{j-1}^{(n)}}}{2}
(Z_{\tau_j^{(n)}\wedge t}-Z_{\tau_{j-1}^{(n)}\wedge t})=\int_0^t g_{s-}\circ dZ_s
\end{equation*}
in the u.c.p.\ topology (Theorem V.26 in Protter \cite{Protter-04}).

The Stratonovich SDE we are interested in is then written in matrix form as 
\begin{equation}
\label{eq:s}
X_t^\circ=x-\int_0^t \nabla U(X_s^\circ)\, ds+\varepsilon \int_0^t F(X_{s-}^\circ) \circ dZ_s,
\end{equation}
or in the coordinate form
\begin{equation*}
(X^\circ)^i_t=x_i-\int_0^t  \partial_{i} U(X_s^\circ)\, ds+
\varepsilon \sum_{j=1}^m \int_0^t F_{ij}(X_{s-}^\circ) \circ dZ_s^j,\quad 1\leq i\leq n.
\end{equation*}
It corresponds to the It\^o SDE
\begin{equation}
\label{eq:si}
X_t^\circ=x-\int_0^t \nabla U(X_s^\circ)\, ds+\varepsilon \int_0^t F(X_{s-}^\circ) 
dZ_s +\frac{\varepsilon^2}{2} \int_0^t F'(X^\circ_{s-}) F(X^\circ_{s-})\, 
d[Z, Z]_s^c 
\end{equation}
which in turn should be understood as
\begin{equation*}
\begin{aligned}
(X^\circ)^i_t=x_i&-\int_0^t  \partial_{i} U(X_s^\circ)\, ds+
\varepsilon \sum_{j=1}^m \int_0^t F_{ij}(X_{s-}^\circ) dZ_s^j \\
&+
\frac{\varepsilon^2}{2} \sum_{l=1}^n\sum_{j,k=1}^m\int_0^t \partial_{l} F_{ij}(X^\circ_{s-}) F_{kl}(X^\circ_{s-})\, 
d[Z^j, Z^k]_s^c,\quad 1\leq i\leq n.
\end{aligned}
\end{equation*} 
Again, from the theory of It\^o SDEs we can conclude that the equation \eqref{eq:s} also has a
 strong solution which is a c\`adl\`ag semimartingale and a strong Markov process (Theorem V.22 in Protter \cite{Protter-04}).

Stratonovich integrals w.r.t.\ martingales are generally not martingales. However, when the integrator is a 
continuous semimartingale, in our case
when $Z=W+\mu t$,
the Stratonovich integral enjoys especially nice properties. 
In this case the solution of  \eqref{eq:s} can be obtained with help of the so-called 
Wong--Zakai approximations. Consider the polygonal approximation $Z^{(n)}$ of the continuous process $Z$ 
\begin{equation}
\label{eq:Zn}
Z^{(n)}_t=Z_{\frac{k}{n}}
+\frac{t-\frac{k}{n}}{n}\Big(Z_{\frac{k+1}{n}}-Z_{\frac{k}{n}}\Big),\quad \frac{k}{n}\leq t\leq \frac{k+1}{n},\ k\geq 0,
\end{equation}
and a path-wise ordinary differential equation
\begin{equation} 
\label{eq:swz}
X_t^{(n)}=x-\int_0^t \nabla U(X_s^{(n)})\, ds+\varepsilon \int_0^t F(X_{s}^{(n)}) \, d Z_s^{(n)}
\end{equation}
with the last integral understood in the Lebesgue--Stieltjes sense.
Then in the limit $n \to\infty$, $X^{(n)}$ converges to the solution $X^\circ$ of \eqref{eq:s},
see Twardowska \cite{Twardowska-96} for a review on the subject and the collection of results.  

Since the processes $Z^{(n)}$ are continuous, the approximation 
 \eqref{eq:swz} and thus the limiting Stratonovich SDE \eqref{eq:s} 
are often chosen for the description of the real world physical processes.

Another remarkable feature of the Stratonovich integral consists in a more simple 
change of variables formula. Let $f\in C^2(\mathbb{R}^n, \mathbb{R})$
and $X^\circ$ be an $n$-dimensional semimartingale. Then (Theorem V.21, \cite{Protter-04})
\begin{equation*}
\begin{aligned}
f(X_t^\circ)=f(X_0^\circ)&+\sum_{j=1}^n\int_0^t \partial_{j} f(X^\circ_{s-})\circ d(X^\circ_s)^j\\
&+\sum_{s\leq t}\Big(f(X_s^\circ)-f(X_{s-}^\circ)-\sum_{j=1}^n \partial_ j f (X_{s-}^\circ)\Delta (X_s^\circ)^j  \Big),
\end{aligned}
\end{equation*} 
so that if the pure jump part of $X^\circ$ vanishes, $(X^\circ)^d\equiv 0$, one obtains the Newton--Leibniz chain rule of Stratonovich integrals. In this case,
 one can construct SDEs on manifolds.

\subsubsection{Canonical Marcus SDE}

Canonical SDEs were introduced by S.~I.~Marcus in \cite{Marcus-78,Marcus-81} in order to preserve the flow property
and a conventional Newton--Leibniz rule for the solutions of SDEs driven by semimartingales with jumps. 
We start with the formulation of the Marcus canonical equation.

The matrix valued function $F$ given,  for any $z\in \mathbb{R}^m$ we
consider the ordinary dif\-fe\-ren\-tial equation
\begin{equation} 
\label{eq:dd} 
\begin{cases}
&\displaystyle \frac{d}{du}y(u)= F(y(u))z,\ u\geq 0,\\
&y(0)=x\in \mathbb{R}^n.
\end{cases}
\end{equation}
Since all $F_{ij}$ are Lipschitz continuous, the vector field $F(\cdot)z$ is complete, that is the solution 
of \eqref{eq:dd} exists and is unique for all $x\in\mathbb{R}^n$ 
and $u\geq 0$. Since $F_{ij}\in C^1(\mathbb{R}^n,\mathbb{R})$, the vector field 
$F(\cdot)z$ generates the flow of diffeomorphisms
\begin{equation*}
\varphi^z_u(x)= y(u, x; z),\quad u\geq 0.
\end{equation*}
We denote $\varphi^z(x):=\varphi^z_1(x)$.

The canonical Marcus SDE is then formally written as
\begin{equation}
\label{eq:wq}
X^{\diamond}(t)=x-\int_0^t \nabla U(X^\diamond_s)\, ds+\varepsilon  \int_0^t F(X_{s-}^\diamond )\diamond d Z_s
\end{equation}
where $\diamond\, dZ$ denotes the Marcus canonical integral. This equation is understood in the following sense:
\begin{equation*}
\begin{aligned}
X^\diamond(t)&=x-\int_0^t \nabla U(X^\diamond_s)\, ds+\varepsilon \int_0^t F(X^\diamond_{s-}) \diamond dZ_s\\
&=x-\int_0^t \nabla U(X^\diamond_s)\, ds+\varepsilon \int_0^t F(X^\diamond_{s-})\circ  dZ^c_s
+\varepsilon \int_0^t F(X^\diamond_{s-})\,  dZ^d_s \\
&+
\sum_{s\leq t} \Big(  \varphi^{\varepsilon \Delta Z_s}(X^\diamond_{s-}) -X^\diamond_{s-}
- F(X^\diamond_{s-}) \varepsilon \Delta Z_s \Big) \\
&=x-\int_0^t \nabla U(X^\diamond_s)\, ds+\varepsilon \int_0^t F (X^\diamond_{s-}) \,dZ_s
+\frac{\varepsilon^2}{2} \int_0^t F'(X^\diamond_{s-})F(X^\diamond_{s-})\,  d[Z,Z]^c_s \\
&+
\sum_{s\leq t} \Big(  \varphi^{\varepsilon \Delta Z_s}(X^\diamond_{s-}) -X^\diamond_{s-}
- F(X^\diamond_{s-})\varepsilon \Delta Z_s\Big),
\end{aligned}
\end{equation*}
where the formula after the second equality sign represents the canonical equation in terms of the Stratonovich integral,
whereas the  formula after the third equality sign gives the It\^o interpretation.

Under above conditions, the canonical equation has
a unique global solution which is a c\`adl\`ag semimartingale, see Theorem 3.2 in Kurtz \textit{et al.} \cite{KurtzPP-95}.
Moreover, this solution is strong Markov (Theorem 5.1 in \cite{KurtzPP-95}).

The jumps of $X^\diamond$ occur only when the jumps of $Z$ occur. If
$Z_s$ does not have a jump at $s$, $\Delta Z_s=0$, then the trajectory $X^\diamond$ moves continuously like the
solution of the Stratonovich SDE driven by $Z^c_s$. If $Z_s$ has a jump $\Delta Z_s$ at
time $s$, then the trajectory of the solution jumps from the point $X^\diamond_{s-}$ to 
$\varphi^{\varepsilon \Delta Z_s} (X^\diamond_{s-})$.
That is, it flies from the point $X^\diamond_{s-}$ along the integral curve of the vector field 
$F(X^\diamond_{s-}) \varepsilon \Delta Z_s$
with infinite speed and lands at $\varphi^{\varepsilon \Delta Z_s} (X^\diamond_{s-})$. 
Then the similar movement repeats inductively.
It is clear, that if $Z^d\equiv 0$, then the Marcus SDE coincides with the Stratonovich SDE. 
If the noise is additive, i.e.\ 
$F=\text{const}$, all three SDEs coincide.

It is necessary to note that the Marcus canonical integral w.r.t.\ $Z$ appearing in the equation \eqref{eq:wq} is not a 
proper integral since it can be defined only for a special class of integrands depending on the driving process $Z$ or on the solutions
 $X^\diamond$ of the SDE, namely for processes
of the type $g(Z_{t-})$ or $g(X^\diamond_{t-})$, 
$g$ being a smooth function.
We refer the reader to Chapter 4.4.5 in Applebaum \cite{Applebaum-09} and Definition 4.1 in 
Kurtz \textit{et al.} \cite{KurtzPP-95} for details.

The Wong--Zakai scheme \eqref{eq:swz} can be considered also for jump processes $Z$. 
In this case, the solutions driven by polygonal approximations \eqref{eq:Zn} converge to the solution of the 
canonical equation $X^\diamond$ in the 
sense of weak convergence of finite dimensional distributions, see Corollary on p.\ 329 in Kunita \cite{Kunita-95}.

The change of variables formula for solutions of Marcus SDEs 
has the form of the conventional Newton--Leibniz rule. 
For $f\in C^2(\mathbb{R}^n,\mathbb{R})$ we have (see Proposition 4.2 in \cite{KurtzPP-95})
\begin{equation}
\label{eq:mch}
f(X_t^\diamond )=f(x)
+\sum_{j=1}^n \int_0^t \partial_{x_j} f(X_{s-}^\diamond )\diamond (X^\diamond_s)^j.
\end{equation}
Similar to continuous Stratonovich SDEs, Marcus canonical SDEs can be considered on smooth manifolds, see
\cite{AppKun93,Fujiwara91,KurtzPP-95,Liao04}.

\section{The main result and examples}

\subsection{The first exit time}

Consider the first exit times of the processes $X, X^\circ$ and $X^\diamond$ from the domain $\mathcal{G}$:
\begin{equation*}
\begin{aligned}
&\tau_x( \varepsilon)=\inf\{t\geq 0: X_t(x)\notin \mathcal{G}\},\\
&\tau_x^\circ( \varepsilon)=\inf\{t\geq 0: X^\circ_t(x)\notin \mathcal{G}\},\\
&\tau_x^\diamond( \varepsilon)=\inf\{t\geq 0: X^\diamond_t(x)\notin \mathcal{G}\}.
\end{aligned}
\end{equation*}
The main result of the paper is presented in the following theorem:
\begin{thm}
Define the sets
\begin{equation*}
\label{eq:E}
E=E^\circ:=\{z\in\mathbb{R}^m:\ F(0)z\notin \mathcal{G} \}\quad \text{and}\quad 
E^\diamond:=\{z\in\mathbb{R}^m:\ \varphi^z(0)\notin \mathcal{G}\}
\end{equation*} 
and suppose that $m(E)=m(E^\circ)>0$ and $m(E^\diamond)>0$.
Then for any $u>-1$ and $x\in\mathcal{G}$ the following limits hold:
\begin{equation*}
\begin{aligned}
\lim_{\varepsilon\to 0}\mathbf{E} e^{-u m(E) H(\varepsilon^{-1}) \tau_x(\varepsilon)}
&=\lim_{\varepsilon\to 0}\mathbf{E} e^{-u m(E^\circ) H(\varepsilon^{-1}) \tau_x^\circ(\varepsilon)}\\
&=\lim_{\varepsilon\to 0}\mathbf{E} e^{-u m(E^\diamond) H(\varepsilon^{-1}) \tau_x^\diamond(\varepsilon)}=
\frac{1}{1+u}.
\end{aligned}
\end{equation*} 
Moreover, there is $\gamma>0$ such that this convergence is uniform over all 
$x\in \mathcal{G}$ with $\operatorname{dist} (x,\partial\mathcal{G})\geq \varepsilon^\gamma$.
\end{thm}
In other words, the appropriately normalised first exit times converge in law to the standard exponential distribution;
there is also convergence of all moments.
\begin{exm}[The sets $E=E^\circ$ and $E^\diamond$ are different] 
\emph{Consider a one-dimensional dynamical system $Y$ perturbed by a
bivariate L\'evy process $Z=(Z^1, Z^2)$.
Let $F=(F_1, F_2)$ with
\begin{equation*}
F_1(x)=\frac{1}{(x+1)^2+1}\quad \text{and}\quad F_2(x)=\frac{1}{(x-1)^2+1}.
\end{equation*}
In this case according to \eqref{eq:E}, the set $E=E^\circ$ is a union of two half-plains (Fig.\ \ref{fig:EE} (l.)), 
\begin{equation*}
E=\{z\in \mathbb{R}^2:\ \frac{z_1+z_2}{2}\geq 1 \text{ or }
\frac{z_1+z_2}{2}\leq -1\}.
\end{equation*} 
The set $E^\diamond$ also consists of two halves (Fig.\ \ref{fig:EE} (r.)),
\begin{equation*}
\mathbf{E}^\diamond=\Big\{z\in\mathbb{R}^2:\int_0^1\frac{du}{z_1 F_1(u)+z_2F_2(u)}
\text{ or }
\int_{0}^{-1}\frac{du}{z_1 F_1(u)+z_2F_2(u)}\in(0,1]\Big\}.
\end{equation*}
\begin{figure}[t]
\centerline{\hbox{ 
\psfig{file=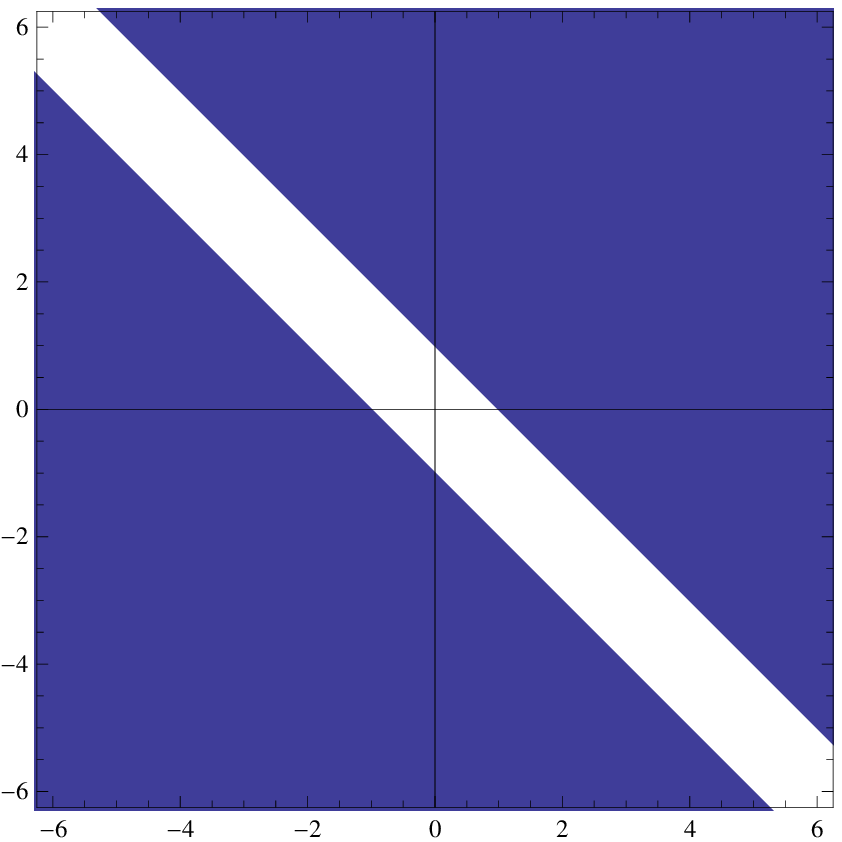,width=0.4\columnwidth} \hspace{0.075\columnwidth}
\psfig{file=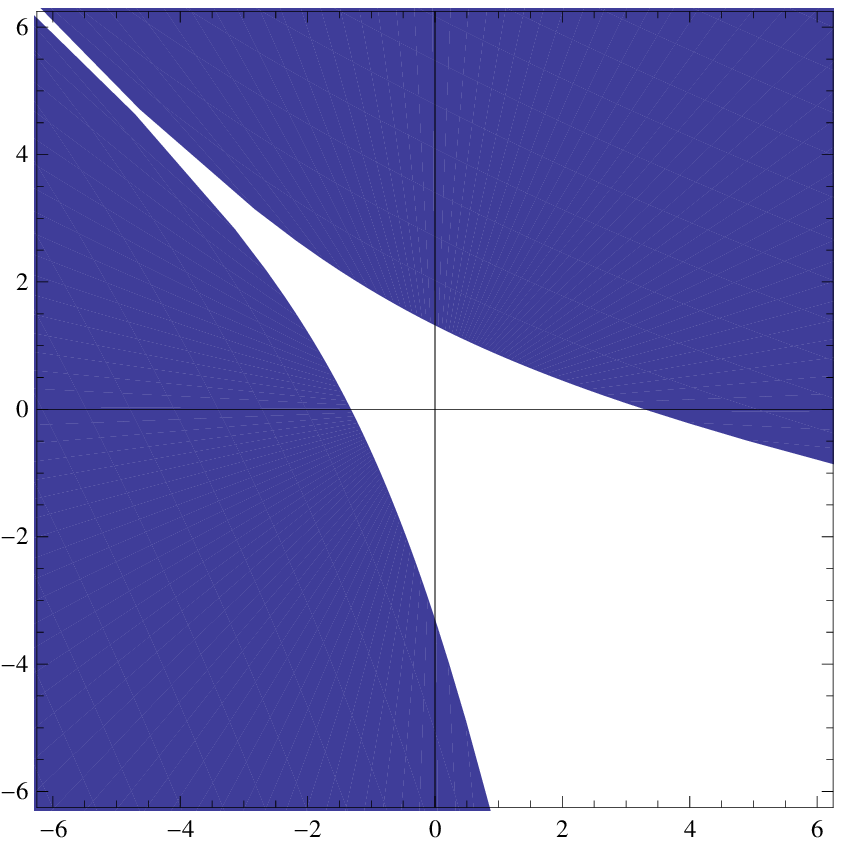,width=0.4\columnwidth}
}}
\caption{The sets $E=E^\circ$ and $E^\diamond$. \label{fig:EE}}
\end{figure}
For example, for a bivariate
isometric Cauchy process $Z$ with the jump measure $\nu(dz)=\|z\|^{-3}\, dz$,
$z\neq 0$, we obtain $H(\varepsilon^{-1})=2\pi\varepsilon$, 
$m(E)\approx 0.49$,
$m(E^\diamond) \approx 0.45$ and hence $\mathbf{E} \tau_x(\varepsilon)\approx \mathbf{E} \tau_x^\circ (\varepsilon)\approx 0.33/\varepsilon$, $ \mathbf{E} \tau_x^\diamond (\varepsilon)
\approx 0.35/\varepsilon$.}
\end{exm}

\begin{exm}[Reduction to additive noise for $n=m=1$]\emph{
In the case $n=m=1$, the exit time $\tau^\diamond$ can be obtained with help of the change of variables formula 
\eqref{eq:mch} and a trick which was used by Nourdin and Simon in \cite{NourdinS-06}.
Consider a one-dimensional canonical Marcus equation 
\begin{equation*}
 X_t^\diamond =X_0-\int_0^t U'(X_{s}^\diamond )\, ds+\varepsilon\int_0^t F(X_{s-}^\diamond) \diamond d Z_s
\end{equation*}
driven by a univariate L\'evy process $Z$
and let $\mathcal{G}=(-a, b)$, $a,b>0$.
Assume that the perturbation is \textit{uniformly elliptic} in $\mathcal{G}$, that is $F(x)>0$, $x\in [-a,b]$, and 
introduce the function 
\begin{equation*}
f(x)=\int_0^x\frac{dy}{F(y)},\quad f'(x)=\frac{1}{F(x)},\quad x\in[-a,b].
\end{equation*}
Applying the change of variables formula \eqref{eq:mch} yields the following SDE for the process $Y^\diamond_t=f(X_t^\diamond)$:
\begin{equation*}
\begin{aligned}
Y^\diamond_t=f(X_t^\diamond)&= f(x)-\int_0^t f'(X_{s}^\diamond)U'(X_{s}^\diamond )\, ds
+\varepsilon\int_0^t f'(X^\diamond_{s-})F(X_{s-}^\diamond) \diamond d Z_s\\
&=f(x)-\int_0^t f'(X_{s}^\diamond)U'(X_{s}^\diamond )\, ds+\varepsilon Z_t\\
&=y^\diamond+\int_0^t B^\diamond (Y_{s}^\diamond )\, ds+\varepsilon Z_t
\end{aligned}
\end{equation*}
with $B^\diamond(y):=-\frac{U'}{F}\circ f^{-1}(y)$ and $y^\diamond=f(x)$. Note that since $F$ is strictly positive and $f$ is monotone 
increasing with $f(0)=0$, we can introduce the new effective potential 
$U^\diamond(y)=-\int_0^y B^\diamond(v)\, dv $ which is a one-well potential with the global minimum at the origin.  
Thus the process $Y^\diamond_t$ satisfies the SDE with \textit{additive} small noise which has been studied in
Imkeller and Pavlyukevich \cite{ImkellerP-06,ImkellerP-08} and Imkeller \textit{et al.} \cite{ImkPavlWetz-09}. 
It is clear that
$X_t^\diamond\notin (-a,b)$ if and only if
$Y^\diamond_t\notin (f^{-1}(a), f^{-1}(b))$.}

\emph{For instance, if $Z$ is a symmetric $\alpha$-stable L\'evy process with $\alpha\in (0,2)$ and the jump measure 
$\nu(dz)=|z|^{-1-\alpha}\,dz$, $z\neq 0$, then for the first exit time of $X^\diamond$ from $(-a,b)$ we immediately
obtain the asymptotics
\begin{equation*}
\lim_{\varepsilon\downarrow 0}\mathbf{E} e^{-u M^\diamond \frac{2}{\alpha}\varepsilon^\alpha \tau_x^\diamond}
=\frac{1}{1+u},\quad u>-1,
\end{equation*}
with
\begin{equation*}
M^\diamond=\Big( \frac{1}{|f^{-1}(-a)|^\alpha}+     \frac{1}{|f^{-1}(b)|^\alpha}   \Big)^{-1}.
\end{equation*}}
\end{exm}

\section{First exit time of the It\^o SDE with multiplicative noise}

\subsection{Big and small jumps of $Z$}

For $\rho\in(0,1)$ and $\varepsilon \leq 1$ let us distinguish the small and big jumps of the driving process $Z$ and decompose it
 into a sum 
\begin{equation*}
Z_t=L_t+\eta_t
\end{equation*}
 with
\begin{equation*}
\eta_t:=\sum_{s\leq t}\Delta Z_s\mathbb{I}\{\|\Delta Z_s\|\geq \varepsilon^{-\rho}\}
\end{equation*}
being a compound Poisson process with the characteristic exponent
\begin{equation*}
\mathbf{E} e^{i\langle \lambda, \eta_t\rangle}
=\exp\Big(   t \int_{\|y\|\geq  \varepsilon^{-\rho}} (e^{i\langle\lambda, y\rangle}-1)\nu(dy)       \Big).
\end{equation*}
The L\'evy process $L$ is a process with bounded jumps, $\|\Delta L_s\|\leq \varepsilon^{-\rho}$, and thus possesses all
moments. Moreover, it is a sum of its continuous component $L^c_t=W_t+\mu t$ being the Brownian motion with drift, 
and a pure jump part $L^d$.

Denote by
$0=\tau_0<\tau_1<\tau_2<\dots$
the successive jump times of $\eta$ and by $J_k$ the respective jump sizes.
The inter-jump times $T_k=\tau_k-\tau_{k-1}$ are iid exponentially distributed random variables with the mean value
\begin{equation*}
\mathbf{E} T_k=\frac{1}{\beta_\varepsilon}:=\Big( \int_{\|y\| \geq \varepsilon^{-\rho}} \nu(dy)\Big)^{\!-1} =\frac{1}{H(\varepsilon^{-\rho})} \to\infty,\quad 
\varepsilon\to 0,
\end{equation*}
and the probability distribution function $\mathbf{P}(T_k\leq  u)=1-e^{-u\beta_\varepsilon}$, $u\geq 0$.
The probability law of $J_k$ is also known explicitly in terms of the L\'evy measure $\nu$:
\begin{equation}
\label{eq:law}
\mathbf{P}(J_k\in A)=\beta_\varepsilon^{-1}\nu(A\cap\{z:\ \|z\|\geq \varepsilon^{-\rho} \}),\quad A\in\mathcal{B}( \mathbb{R}^m).
\end{equation}

\subsection{Perturbations by the process $\varepsilon L$\label{ss:pert}}
 
\begin{lem}
\label{l:q}
Let $\rho\in (0,1)$, $\mu_\varepsilon:=\mathbf{E} L_1$ and $T_\varepsilon:=\varepsilon^{-\theta}$ for some $\theta>0$.
There exist $\delta_0=\delta_0(\rho)>0$ and $\theta_0=\theta_0(\rho)>0$ such that for all 
$\delta\in (0,\delta_0)$, $\theta\in (0,\theta_0)$ there are $p_0=p_0(\delta)$ and 
$\varepsilon_0=\varepsilon_0(\rho,\delta)$ such that the estimates
\begin{equation*}
 \|\varepsilon L_{T_\varepsilon}\|=\varepsilon \|\mu_\varepsilon\| T_\varepsilon <   \varepsilon^{2\delta}\quad \text{and}\quad 
\mathbf{P}([\varepsilon L]_{T_\varepsilon}^d\geq \varepsilon^\delta)\leq \exp(-\varepsilon^{-p})
\end{equation*}
 hold for all $p\in (0,p_0)$ and $0<\varepsilon\leq \varepsilon_0$. 
\end{lem}
\begin{proof}
Let us represent the process $L$ as 
\begin{equation*}
L_t=\tilde L_t+\mu_\varepsilon t
\end{equation*}
with $\tilde L$ being a zero mean L\'evy martingale with bounded jumps.

\noindent
\textbf{Step 1.} We have the following estimates for the mean value $\mu_\varepsilon$:
\begin{equation*}
\begin{aligned}
&\mu_\varepsilon^i :=\mathbf{E} L^i_1=\int_{1<\|z\|\leq \varepsilon^{-\rho}} z_i\,\nu(dz),\quad 1\leq i\leq m,\\
&\|\mu_\varepsilon\|^2\leq \int_{1<\|z\|\leq \varepsilon^{-\rho}} \|z\|^2 \,\nu(dz)=-\int_{1}^{\varepsilon^{-\rho}} u^2 \, dH(u)
\leq \varepsilon^{-2\rho}H(1),\\
&\|\mu_\varepsilon\|\leq \sqrt{H(1)}\varepsilon^{-\rho}.  
\end{aligned}
\end{equation*}
Consequently, for any $\rho\in (0,1)$, $\theta_0:=(1-\rho)/3$ and $\delta_0:=(1-\rho)/4$
we obtain
\begin{equation*}
\varepsilon\|\mu_\varepsilon\|T_\varepsilon< \varepsilon^{2\delta}
\end{equation*}
for all $0<\delta<\delta_0$ and $0<\theta<\theta_0$ and $\varepsilon$ small enough.

\noindent 
\textbf{Step 2.}
The of quadratic variation process $[\varepsilon \tilde L]^d_t$  is a L\'evy subordinator
\begin{equation*}
{}[\varepsilon \tilde L]^d_t=\varepsilon^2\sum_{s\leq t} \|\Delta \tilde L\|_s^2 
=\varepsilon^2\int_0^t \int_{0<\|z\|\leq \varepsilon^{-\rho}} \|z\|^2 N(dz, ds).
\end{equation*}
Since the jumps of $[\varepsilon \tilde L]^d$ are bounded,
its Laplace transform is well-defined for all $\lambda\in \mathbb{R}$ and equals 
\begin{equation}
\label{eq:22} 
\begin{aligned}
\mathbf{E} e^{\lambda [\varepsilon \tilde L]^d_t}&=
\exp\Big( t\int_{0<\|z\|\leq \varepsilon^{-\rho}} (e^{\lambda \varepsilon^2\|z\|^2}-1)\, \nu(dz)\Big)\\
&=\exp\Big(- t\int_{0<u\leq \varepsilon^{-\rho}} (e^{\lambda \varepsilon^2 u^2}-1)\, dH(u)\Big).
\end{aligned}
\end{equation}
For any $\lambda>0$, the exponential Chebyshev inequality implies  that
\begin{equation}
\begin{aligned} \label{eq:df}
\mathbf{P}([\varepsilon \tilde L]^d_{T_\varepsilon}> \varepsilon^{\delta})&
=\mathbf{P}(e^{ \lambda[\varepsilon \tilde L]^d_{T_\varepsilon}}> e^{\lambda \varepsilon^{\delta}})
\leq e^{-\lambda\varepsilon^{\delta}}\mathbf{E}  e^{\lambda[\varepsilon \tilde L]^d_{T_\varepsilon}}\\
&=\exp\Big(-\lambda\varepsilon^{\delta}- T_\varepsilon\int_{0<u\leq \varepsilon^{-\rho}} (e^{\lambda \varepsilon^2 u^2}-1)\, dH(u)\Big).
\end{aligned}
\end{equation}
For $\lambda=\lambda_\varepsilon:=\varepsilon^{-2\delta}$ with $0<\delta<\delta_0=(1-\rho)/4$ we have  
$\max_{0\leq u\leq \varepsilon^{-\rho}} \lambda_\varepsilon \varepsilon^{2} u^2\leq \lambda_\varepsilon \varepsilon^{2(1-\rho)}\downarrow 0$ as 
$\varepsilon \downarrow 0$. With help of the elementary inequality $e^x-1\leq 2x$ for small positive $x$
the second summand  appearing in the exponent in r.h.s.\ of \eqref{eq:df} can be now estimated as
\begin{equation*}
\begin{aligned}
\Big|T_\varepsilon\int_{0<u\leq \varepsilon^{-\rho}}& (e^{ \lambda_\varepsilon \varepsilon^2 u^2}-1)\, dH(u)\Big|
\leq \Big|  2T_\varepsilon \lambda_\varepsilon \varepsilon^2  \Big(\int_{0<u\leq 1} 
+\int_{1<u\leq \varepsilon^{-\rho}}\Big)u^2 \, dH(u)\Big|\\
&\leq   2T_\varepsilon \lambda_\varepsilon \varepsilon^2 \Big|\int_{0<u\leq 1}u^2\, dH(u) \Big|
+ 2T_\varepsilon \lambda_\varepsilon \varepsilon^{2(1-\rho)}\Big|\int_{1<u\leq \varepsilon^{-\rho}} \, dH(u)\Big|\\
&\leq  C T_\varepsilon \lambda_\varepsilon \varepsilon^2  
+2 H(1) T_\varepsilon \lambda_\varepsilon \varepsilon^{2(1-\rho)}
\end{aligned}
\end{equation*}
with $C=\big|\int_{0<u\leq 1}u^2\, dH(u)\big|>0$.
Consequently, for all $0<\delta<\delta_0$ and $0<\theta<\theta_0$ we see that
the exponential inequality
\begin{equation*}
\begin{aligned}
\mathbf{P}([\varepsilon \tilde L]^d_{T_\varepsilon}> \varepsilon^{\delta})&\leq 
\exp\Big(-\lambda_\varepsilon\varepsilon^\delta+C T_\varepsilon\lambda_\varepsilon \varepsilon^2  
+2H(1) T_\varepsilon\lambda_\varepsilon \varepsilon^{2(1-\rho)}\Big)
\leq e^{-\varepsilon^{-\delta/2}}
\end{aligned}
\end{equation*}
holds for $\varepsilon$ small enough and the Lemma holds with $p\in(0,\delta/2)$.
\end{proof} 

\begin{lem}
\label{l:qq}
Let $\rho\in (0,1)$ and $(g_t)_{t\geq 0}$ be a bounded adapted c\`adl\`ag stochastic process mit values in $\mathbb{R}^m$,
$T_\varepsilon=\varepsilon^{-\theta}$, $\theta>0$.
There are $\delta_0=\delta_0(\rho)>0$ and $\theta_0=\theta_0(\rho)>0$ such that for all $\delta\in (0,\delta_0)$ and 
$\theta\in (0,\theta_0)$ there are $p_0=p_0(\rho,\delta)$ and $\varepsilon_0=\varepsilon_0(\delta)$ such that
the exponential estimate
\begin{equation*}
\mathbf{P}\Big(\sup_{0\leq t\leq T_\varepsilon}\varepsilon
\Big|\sum_{j=1}^m\int_0^t g^j_{s-}\, d\tilde L^j_s  \Big|\geq \varepsilon^\delta \Big)\leq e^{-\varepsilon^{-p}}.
\end{equation*}  
holds for all $p\in (0,p_0)$ and $0<\varepsilon\leq \varepsilon_0$.
\end{lem}

\begin{proof}
\textbf{Step 1.} Suppose that $\sup_{t\geq 0}\|g_t\|\leq C$ for some $C>0$.
Consider the one-dimensional martingale
\begin{equation*}
M_t=\sum_{j=1}^m\int_0^t g^j_{s-}\, d\tilde L^j_s.
\end{equation*}
By construction $|\Delta M_t|\leq C\varepsilon^{-\rho}$.
We estimate the probability of a deviation of the size $\varepsilon^\delta$ of $\varepsilon M_t$ from zero with help of the 
exponential inequality for martingales, 
see Theorem 26.17(i) in Kallenberg \cite{Kallenberg-02}. Indeed for any $\delta>0$ and $\theta>0$ we have
\begin{equation*}
\mathbf{P}\Big(\sup_{t\leq T_\varepsilon}| \varepsilon M_t | \geq \varepsilon^{\delta}\Big)\leq 
\mathbf{P}\Big(\sup_{t\leq T_\varepsilon}| \varepsilon M_t | \geq \varepsilon^{\delta}\Big| [\varepsilon M]_{T_\varepsilon}\leq \varepsilon^{4\delta} \Big)
+\mathbf{P}([\varepsilon M]_{T_\varepsilon}> \varepsilon^{4\delta}).
\end{equation*}
Inspecting the proofs of Lemma 26.19 and Theorem 26.17(i) in Kallenberg \cite{Kallenberg-02} we get that for any
$\lambda>0$ 
\begin{equation*}
\mathbf{P}\Big(\sup_{t\leq T_\varepsilon}\varepsilon M_t\geq  \varepsilon^{\delta}\Big|[\varepsilon M]_{T_\varepsilon}\leq \varepsilon^{4\delta} \Big)
\leq e^{- \lambda \varepsilon^{\delta}+\lambda^2 h( \lambda C \varepsilon^{1-\rho})\varepsilon^{4\delta}}
\end{equation*}
with $h(x)=-(x+\ln(1-x)_+)x^{-2}$. 
For any $0<\delta<\delta_1:=(1-\rho)/2$ we set $\lambda=\lambda_\varepsilon=\varepsilon^{-2\delta}$, so that 
$h( \lambda_\varepsilon C \varepsilon^{1-\rho})\to 1/2$ as $\varepsilon\to 0$. Hence
we obtain the estimate
\begin{equation*}
\mathbf{P}\Big(\sup_{t\leq T_\varepsilon}\varepsilon M_t\geq  \varepsilon^{\delta}\Big|[\varepsilon M]_{T_\varepsilon}
\leq \varepsilon^{4\delta} \Big)\leq e^{- \varepsilon^{-\delta/2}}\leq e^{- \varepsilon^{-p}},
\end{equation*}
which holds for $\varepsilon$ small enough and $p\in(0,\delta/2)$.
An analogous inequality holds for $\inf_{t\leq T_\varepsilon}\varepsilon M_t$.

\noindent
\textbf{Step 2.}
There is a constant $C_1>0$ with
\begin{equation*}
{}[\varepsilon M]_t=\int_0^t g_s^2 \, d [\varepsilon W]_s+\int_0^t g_{s-}^2 \, d [\varepsilon\tilde L]^d_s\leq C_1( \varepsilon^2 t+ [\varepsilon \tilde L]^d_t),\quad
t\geq 0.
\end{equation*}
For sufficiently small $\varepsilon$ we obtain the estimate
\begin{equation*}
\mathbf{P}([\varepsilon M]_{T_\varepsilon} \geq \varepsilon^{4\delta})\leq \mathbf{P}([\varepsilon\tilde L]^d_{T_\varepsilon} \geq \varepsilon^{5\delta})
+ \mathbf{P}(C_1\varepsilon^2 T_\varepsilon\geq\varepsilon^{5\delta} ). 
\end{equation*}
For $\varepsilon$ small enough the second summand vanishes for all $0<\delta<\delta_2$, $0<\theta<\theta_2$ with 
$0<\theta_2+5\delta_2<2$. The first summand is bounded by $e^{-\varepsilon^{-p}}$ from above due
to Lemma \ref{l:q} for  $0<\delta<\delta_3$ some $0<\theta<\theta_3$ and $0<p<p_1(\rho,\delta)$ and $\varepsilon$ 
small enough. This the statement of the Lemma holds with $\delta_0=\min\{\delta_1,\delta_2,\delta_3\}$, 
$\theta_0=\min\{\theta_1,\theta_2,\theta_3\}$ and $p_0=\min\{\delta/2, p_1(\delta)\}$.
\end{proof}

\begin{lem}
\label{l:p}
For any $\rho\in (0,1)$ there are  $\gamma_0>0$ and  $p>0$ such that for all $0<\gamma\leq \gamma_0$
\begin{equation*}
\sup_{x\in \mathcal{G}} \mathbf{P}\Big(\sup_{0\leq  t< T_1}\| X(t,x)-Y(t,x)\|\geq  \varepsilon^\gamma \Big) \leq e^{-\varepsilon^{-p}}.
\end{equation*} 
\end{lem}

\begin{proof}
\textbf{Step 1.} By assumptions on the potential $U$, any deterministic trajectory $Y_t(y)$, $y\in \mathcal{G}$, reaches 
a small fixed neighbourhood of the origin in some finite time. After entering this small neighbourhood, the trajectory 
$Y_t(y)$ is attracted to the origin with the speed approximately proportional to $C_1\|Y_t\|$, 
$C_1>0$ being the smallest eigenvalue of the matrix $\frac{\partial^2 }{\partial x_i\partial x_j}U(x)\big|_{x=0}$. 
This allows us to estimate of increase rate of the time a trajectory $Y_t(y)$, $y\in \mathcal{G}$, needs to reach some 
$\varepsilon^\delta$-neighbourhood of the origin. Indeed, for any $\delta>0$ the following inequality holds true
for any $0\leq \varepsilon\leq \varepsilon_0$, $y\in\mathcal{G}$ and $\varepsilon_0>0$ small enough: 
\begin{equation}
\label{eq:99}
\|Y_t(y)\|\leq \varepsilon^\delta,\quad t\geq V_\varepsilon:=\frac{2\delta}{C_1}|\ln \varepsilon|.
\end{equation} 
\textbf{Step 2.} 
Here we show that on the time intervals up to $V_\varepsilon$, the random trajectory $X^\varepsilon(x)$ does not deviate much 
from the deterministic solution $Y_t(x)$
with the same initial value in the absence of big jumps of the driving process $Z$.

For $x\in \mathcal{G}$, with help of Gronwall's lemma we estimate
\begin{equation*}
\begin{aligned}
&\|X_{t\wedge V_\varepsilon\wedge  T_1-}-Y_{t\wedge V_\varepsilon\wedge  T_1-}\|\\
&\qquad\qquad\leq  C_{\text{Lip}}\int_0^{t\wedge V_\varepsilon\wedge  T_1-} \|X_s-Y_s\|\, ds
+ \varepsilon\Big\| \int_0^{t\wedge V_\varepsilon\wedge  T_1-} F(X_{s-})\, dZ_s  \Big\|, \\
&\sup_{0\leq t< V_\varepsilon\wedge T_1}\|X_t-Y_t\|\leq e^{C_{\text{Lip}} V_\varepsilon}
\sup_{0\leq t\leq  V_\varepsilon} 
 \varepsilon\Big\| \int_0^{t} F(X_{s-})\, dL_s  \Big\|.
\end{aligned}
\end{equation*}
Recalling the definition of $V_\varepsilon$ in \eqref{eq:99} and taking into account that 
$V_\varepsilon\leq T_\varepsilon=\varepsilon^{-\theta}$ for 
$\varepsilon$ small enough and any $\theta>0$
we get for any $\delta>0$ that
\begin{equation*}
\begin{aligned}
\mathbf{P}\Big(\sup_{0\leq t < V_\varepsilon\wedge T_1}&\|X_t-Y_t\|\geq \varepsilon^\delta\Big)
\leq \mathbf{P}\Big(
\sup_{0\leq t\leq  V_\varepsilon} \varepsilon\Big\| \int_0^t F(X_{s-})\, dL_s  \Big\|\geq \varepsilon^{(1+\frac{2C_{\text{Lip}}}{C_1})\delta}\Big)\\
&\leq \mathbf{P}\Big(
\sup_{0\leq t\leq  T_\varepsilon} \varepsilon\Big\| \int_0^t F(X_{s-})\, d\tilde L_s  \Big\|+C_2\varepsilon T_\varepsilon\|\mu_\varepsilon\| 
\geq \varepsilon^{(1+\frac{2C_{\text{Lip}}}{C_1})\delta}\Big)
\end{aligned}
\end{equation*}
with some $C_2>0$. With help of 
Lemmas \ref{l:q} and \ref{l:qq} we find that there are $\delta_1>0$ and $\theta_1>0$ such that the 
last probability is smaller than $\exp(-\varepsilon^{-p})$ for all $\delta\in (0,\delta_1)$,  $\theta\in(0,\theta_1)$ and $p\in(0,p_1(\delta))$.

\smallskip
\noindent
\textbf{Step 3.}
In this Step we exploit the attractor property of the origin and show that in the absence of big jumps of the driving 
process $Z$ the random path $X_t(x)$ with the initial value close to the origin does not deviate much
on the polynomially long time intervals $T_\varepsilon$.

Consider the function $f(x)=\ln(1+U(x))\geq 0$. For $\|x\|$ small, one can estimate  $c_1 \|x\|^2\leq f(x)\leq c_2\|x\|^2$
for some positive $c_1$ and $c_2$. 
Furthermore, the derivatives $\partial_i f(x)=\frac{\partial_i U(x)}{1+U(x)}$,
$\partial_i\partial_j f(x)=2\frac{\partial_i U(x)\partial_j U(x)-\partial_i\partial_j U(x)}{(1+U(x))^2}$ are bounded due to 
assumptions on $U$.

We apply the It\^o formula to the process $f(X_t)$:
\begin{equation*}
\begin{aligned}
0&\leq f(X_{t\wedge T_\varepsilon\wedge T_1-})\\
&\begin{aligned}
=f(x) &+\sum_{i=1}^n \int_0^{t\wedge T_\varepsilon\wedge T_1-}\!\! \partial_i f(X_{s-})\, dX^i_s 
+\frac{1}{2}\sum_{i,j=1}^n\int_0^{t\wedge T_\varepsilon\wedge T_1-}\!\! \partial_i\partial_j f(X_{s-})\, d [X^i,X^j]^c_s\\
&+\sum_{s< {t\wedge T_\varepsilon\wedge T_1}} \Big(f(X_{s})-  f(X_{s-}) -\sum_{i=1}^n\partial_if(X_{s-})\Delta X^i_s  \Big)\\
=f(x) &-\int_0^{t\wedge T_\varepsilon\wedge T_1-}\!\frac{\|\nabla U(X_{s-})\|^2}{1+U(X_{s-})}\, ds 
 +\varepsilon \sum_{i,j=1}^{n,m}\!\int_0^{t\wedge T_\varepsilon\wedge T_1-} 
 \frac{\partial_i U(X_{s-}) F_{ij}(X_{s-})}{1+U(X_{s-})}\, dZ^j_s  \\
&+\frac{\varepsilon^2}{2}\sum_{i,j=1}^n \sum_{k,l=1}^m \int_0^{t\wedge T_\varepsilon\wedge T_1-} 
 \partial_i\partial_j f(X_{s-})F_{ik}(X_{s-})F_{jl}(X_{s-})\, d [Z^k, Z^l]^c_s\\
&+\sum_{s< {t\wedge T_\varepsilon\wedge T_1}} \Big(  f(X_{s})-  f(X_{s-})   -\sum_{i=1}^n\partial_if(X_{s-})\Delta X^i_s  \Big).
\end{aligned}
\end{aligned}
\end{equation*}
We note that the first integral in the last formula is non-negative, the integrands in the It\^o integral w.r.t.\ 
$Z$ and in the integrals w.r.t.\ 
$[Z^k, Z^l]^c$ are bounded, the quadratic covariations satisfy
 $[Z^k, Z^l]^c_t=[W^k, W^l]_t=\sigma_{kl}t$, and finally the estimate
\begin{equation*}
\begin{aligned}
\sum_{s\leq t}\Big| f(X_{s})&-f(X_{s-})-\sum_{i=1}^n  \partial_i f(X_{s-})\Delta X^i_s \Big|\\
&\leq \frac{1}{2}\sum_{i,j=1}^n \sum_{s\leq t}\Big| \int_0^1  (1-v) \partial_i\partial_j f(X_{s-}+v\Delta X_s) \, dv\Big|
\cdot|\Delta X^i_s\Delta X^j_s|\\
&\leq C_3 \sum_{s\leq t} \|\Delta X_s\|^2=C_3 [X]_t^d
\end{aligned}
\end{equation*}
holds
with some $C_3>0$.
Furthermore, since $F$ is bounded the estimate
\begin{equation*}
{}[X]_t^d\leq C_4[Z]^d_t=C_4[L]^d_t
\end{equation*}
holds
for some constant $C_4>0$ and $0\leq t<T_1$. 
Combining these estimates and denoting $g(x)= \frac{\nabla^T U(x)F(x)}{1+U(x)}$ we obtain the following estimate with 
some positive constant $C_5$, $\varepsilon$ small enough and $\|x\|\leq \varepsilon^\delta$, $\delta>0$:
\begin{equation*}
0\leq f(X_{t\wedge T_\varepsilon\wedge T_1-})
\leq C_5\Big( \varepsilon^{2\delta}  +\varepsilon \sup_{0\leq t\leq T_\varepsilon}\Big|\int_0^t g(X_{s-})\, d\tilde L_s\Big|  
+\varepsilon  \|\mu_\varepsilon\|T_\varepsilon
+  \varepsilon^2 {T_\varepsilon}      +\varepsilon^2 [L]^d_{T_\varepsilon}\Big).
\end{equation*} 
Let the estimates of the Lemmas \ref{l:q} and  \ref{l:qq} hold simultaneously 
for  $\delta\in (0,\delta_2)$, $\theta\in (0,\theta_2)$ and $p\in(0,p_2(\delta))$ for 
some positive $\delta_2$, $\theta_2$, $p_2(\delta)$ and $\varepsilon$ small.  
Then for $\delta\in (0,\delta_2/3)$, $\theta\in (0,\theta_1)$
and $\varepsilon$ small we get
 \begin{equation*}
\begin{aligned}
\mathbf{P}\Big(\sup_{0\leq t< T_\varepsilon\wedge T_1} &\|X_t\|\geq \varepsilon^\delta\Big)\leq
 \mathbf{P}\Big( \varepsilon^{2\delta}\geq \varepsilon^{3\delta}\Big)+
\mathbf{P}\Big( \varepsilon \sup_{0\leq t\leq T_\varepsilon}\Big|\int_0^t g(X_{s-})\, d\tilde L_s\Big|     \geq \varepsilon^{3\delta}\Big)\\
&+ \mathbf{P}\Big( \varepsilon  \|\mu_\varepsilon\|T_\varepsilon   \geq  \varepsilon^{3\delta}\Big)  +\mathbf{P}\Big( \varepsilon^2 {T_\varepsilon}  \geq  \varepsilon^{3\delta}\Big)
+\mathbf{P}\Big(  \varepsilon^2 [L]^d_{T_\varepsilon}   \geq \varepsilon^{3\delta}\Big)\leq e^{-\varepsilon^{-p}}.
\end{aligned}
\end{equation*}
for  $p\in(0,p_2(\delta)/2)$.

\smallskip
\noindent
\textbf{Step 4.}
Combining the estimates of Steps 1, 2 and 3, we  extend the estimate of the Step 3 to all initial values $x\in\mathcal{G}$:
 \begin{equation*}
\mathbf{P}\Big(\sup_{0\leq t< T_\varepsilon\wedge T_1}\|X_t-Y_t\|\geq \varepsilon^\delta\Big)\leq e^{-\varepsilon^{-p}}
\end{equation*}
for $\delta\in (0,\delta_3)$, $\theta\in (0,\theta_3)$, $p\in(0,p_3(\delta))$ and $\varepsilon$ small enough.
Here, $\delta_3=\min\{\delta_1,\delta_2/3\}$,  $\theta_3=\min\{\theta_1,\theta_2\}$ and $p_3=\min\{p_1(\delta), p_2(\delta)/2\}$.

\smallskip
\noindent
\textbf{Step 5.} In this final Step we extend the estimate of the Step 4 from the time interval 
$[0, T_\varepsilon\wedge T_1)$ to the time interval 
$[0, T_1)$.

Denote $X^L$ the solution of the SDE \eqref{eq:i} driven by the process $L$. Clearly, $X_t=X_t^L$ on the event
$\{t<T_1\}$. 
Let $x\in\mathcal{G}$ and $k\geq 1$, then for any $\delta>0$ and $\theta>0$ we have
\begin{equation*}
\begin{aligned}
\mathbf{P}\Big(\sup_{t\in [0, T_1)} \|X_t-Y_t\|\geq \varepsilon^\delta\Big)
&\leq \mathbf{P}\Big(\sup_{t\in [0, kT_\varepsilon \wedge T_1)} \|X_t-Y_t\|\geq \varepsilon^\delta\Big)+
\mathbf{P}(T_1\geq kT_\varepsilon)\\
&\leq \mathbf{P}\Big(\sup_{t\in [0, kT_\varepsilon]} \|X_t^L-Y_t\|\geq \varepsilon^\delta\Big)+
\mathbf{P}(T_1\geq kT_\varepsilon).
\end{aligned}
\end{equation*} 
Moreover for $\varepsilon$ small enough  we have $\|Y(T_\varepsilon,x)\|<\varepsilon^{2\delta}$ for $x\in\mathcal{G}$.
Denote
\begin{equation*}
A_j=\{\sup_{t\in [jT_\varepsilon, (j+1) T_\varepsilon]} \|X_t^L-Y(t-jT_\varepsilon;X_{jT_\varepsilon}^L )\|< \varepsilon^\delta\},\quad 0\leq j\leq k-1.
\end{equation*}
In particular, the probability of $A_0^c=\{\sup_{t\in [0,  T_\varepsilon]} \|X_t^L-Y_t\|\geq  \varepsilon^\delta\}$ was estimated in Step 4. 
Further, for any $k\geq 1$
\begin{equation*}
\bigcap_{j=0}^{k-1} A_j
\subseteq  \Big\{\sup_{t\in [0, kT_\varepsilon]} \|X_t^L-Y_t\|<2 \varepsilon^\delta\Big\}.
\end{equation*}
Consequently
\begin{equation}
\label{eq:33}
\begin{aligned}
\mathbf{P}\Big(\sup_{t\in [0, kT_\varepsilon]} \|X^L_t-Y_t\|&\geq 2\varepsilon^\delta\Big) 
\leq \mathbf{P}\Big(  \bigcup_{j=0}^{k-1}A_j^c    \Big) \\
&=\mathbf{P}\Big(  A_0^c\cup (A_0A_1^c) \cup \cdots\cup (A_0\cdots A_{k-2}A_{k-1}^c ) \Big) \\
&\leq \sum_{j=0}^{k-1}\mathbf{P}(A_j^c, \|X^L_{jT_\varepsilon}\|\in\mathcal{G})\leq k \sup_{x\in\mathcal{G}}\mathbf{P}(A_0^c).
\end{aligned}
\end{equation}
For $k=k_\varepsilon=[\varepsilon^{-2r}]$ and any $\theta>0$ we have
\begin{equation*}
\mathbf{P}(T_1\geq k_\varepsilon T_\varepsilon)=e^{-k_\varepsilon T_\varepsilon \beta_\varepsilon}\leq \exp(-\varepsilon^{r\rho-\theta-2r} l(\varepsilon^{-\rho}))\le e^{-\varepsilon^{-p}}.
\end{equation*}
for all $0<p<p_4:=(2-\rho)r$ and $\varepsilon$ small.
On the other hand,  \eqref{eq:33} and Step 4 yield
\begin{equation*}
\mathbf{P}(\sup_{t\in [0, kT_\varepsilon]} \|X^L_t-Y_t\|\geq 2\varepsilon^\delta) \le \varepsilon^{-2r} e^{-\varepsilon^{-p}}\le  e^{-\varepsilon^{-p/2}}
\end{equation*}
with any $0<p<p_3(\delta)$.
Consequently, the statement of the Lemma holds for any $0<\gamma< \delta_3$, $0<p<\min\{p_3(\gamma)/2,p_4\}$ 
and $\varepsilon$ small enough.
\end{proof}

\subsection{The first exit time of solutions of the It\^o SDE\label{s:Ito}}

Having established the key estimates about deviations of the random  trajectory $X_t(x)$ from the deterministic 
path $Y_t(x)$ on random time intervals between big jumps of the driving process $Z$ we can calculate the asymptotics 
of the Laplace transform of the first exit time. The proof here goes along the lines of the one-dimensional case considered in
Imkeller and Pavlyukevich \cite{ImkellerP-06a} and the multivariate case of a dynamical system driven by a 
multifractal $\alpha$-stable
noise considered in Imkeller \textit{et al.} \cite{ImkPavSta-10}. For the sake of completeness we 
briefly sketch the main idea of the proof.

The argument is based on the concept of the \textit{one big jump} which is often used in the study of 
heavy tail phenomena.
Roughly speaking, it can be shown that under certain conditions the small perturbations of the dynamical system $Y$ 
due the process $\varepsilon L$ can be neglected, and the exit from the domain occurs with high probability at one of the 
jump times $\tau_k$. Just before the time $\tau_k$ the solution $X$
stays in a small neighbourhood of the stable point, so the exit occurs if the jump $\varepsilon J_k$ is large enough, namely if 
$F(X_{\tau_k-})\varepsilon J_k\approx F(0)\varepsilon J_k \notin \mathcal{G}$. The events
$\{\varepsilon J_1\notin E\}=\{F(0)\varepsilon J_1 \in \mathcal{G}\}$,\dots, 
$\{\varepsilon J_{k-1}\notin E\}=\{F(0)\varepsilon J_{k-1} \in \mathcal{G}\}$ and  
$\{\varepsilon J_k\in E\}=\{F(0)\varepsilon J_k \notin \mathcal{G}\}$ 
are independent and build up a geometric sequence of events. 
Their probabilities can be calculated in the limit of $\varepsilon\to 0$ with help of the scaling property of the jump measure $\nu$.

The statement of the main theorem follows from the small noise estimates 
from below and above of the Laplace transforms of the normalised first exit times.
Here we consider the less complicated estimate form below.

For any $u> -1$, with help of the formula of the total probability we have
\begin{equation}
\label{eq:tp}
\mathbf{E} e^{-u m(E) H(\varepsilon^{-1}) \sigma_x(\varepsilon)}\geq \sum_{k=1}^\infty 
\mathbf{E} \Big[e^{-u m(E) H(\varepsilon^{-1}) \tau_k}\mathbb{I}\{\sigma=\tau_k\}\Big].
\end{equation}
For any $\delta>0$ small enough denote $\mathcal{G}^{-\delta}:=\{x\in \mathcal{G}:\ \operatorname{dist}(\partial G,x)\geq \delta\}$ 
the inner part of $\mathcal{G}$ and $\mathcal{G}^{+\delta}\equiv \mathcal{G}^{\delta}:=\{x\in \mathbb{R}^n:\ \operatorname{dist}(G,x)\le \delta\}$ be the outer $\delta$-neighbourhood.
For $k\geq 1$, the strong Markov property allows to write
\begin{equation}
\label{eq:tt}
\begin{aligned}
 \mathbf{E}&\Big[e^{-u m(E) H(\varepsilon^{-1})\tau_k}\mathbb{I}\{\sigma=\tau_k\}\Big]\\
&=\mathbf{E}\Big[ e^{-u m(E) H(\varepsilon^{-1}) \tau_k}\mathbb{I}\{X_t \in \mathcal{G},  t \in [0,\tau_k), X_{\tau_k} \notin \mathcal{G}\}\Big]\\
 &= \mathbf{E}\Big[\prod_{j=1}^{k-1}e^{-u m(E) H(\varepsilon^{-1})T_j}\mathbb{I}\{X_{t+\tau_{j-1}}\in \mathcal{G}, t \in [0,T_{j}]\}\\
&\qquad\quad\times
e^{-u m(E) H(\varepsilon^{-1})T_k}\mathbb{I}\{X_{t+\tau_{k-1}}\in \mathcal{G}, t \in [0,T_k)\}\mathbb{I}\{X_{\tau_k} \notin \mathcal{G}\}\Big]\\
&\geq  
\Big(\inf_{y\in\mathcal{G}^{-\delta}} 
\mathbf{E}\Big[e^{-u m(E) H(\varepsilon^{-1})T_1}\mathbb{I}\{X_{t}(y)\in \mathcal{G}^{-\delta}, t \in [0,T_{1}]\}\Big]\Big)^{k-1}\\
&\qquad\quad\times
\inf_{y\in\mathcal{G}^{-\delta}} \mathbf{E} \Big[
e^{-u m(E) H(\varepsilon^{-1})T_1}\mathbb{I}\{X_{t}(y)\in \mathcal{G}  , t \in [0,T_1)\}\mathbb{I}\{X_{T_1} \notin \mathcal{G}\}\Big].
\end{aligned}
\end{equation}
Let $\gamma>0$ be such that the estimates of the Section \ref{ss:pert} hold. We set $\delta:=\delta(\varepsilon)=\varepsilon^\gamma$.
The exit from the domain with a big jump $\varepsilon J_1$ occurs when $F(X_{T_1-})\varepsilon J_1\notin \mathcal{G}$. 
Further, $\sup_{0\le t< T_1}\|X_{t}-Y_{t}\|\le \frac{1}{2}\varepsilon^\gamma$ with probability exponentially close to 1 
(Lemma \ref{l:p}), 
$Y_t(x)$ reaches a $\frac{1}{2}\varepsilon^\gamma$-neighbourhood of the origin during the relaxation time $V_\varepsilon=\mathcal{O}(|\ln \varepsilon|)$, 
and $T_1> V_\varepsilon$ with high probability. Taking into account that $T_1$ is exponentially distributed with the parameter $\beta_\varepsilon$
we calculate the Laplace transform of $m(E)H(\varepsilon^{-1})T_1$ explicitly, namely
\begin{equation*}
\mathbf{E} e^{-u m(E) H(\varepsilon^{-1})T_1}=\frac{\beta_\varepsilon}{\beta_\varepsilon+ u m(E)H(\varepsilon^{-1})}
=\frac{1}{1+u a_\varepsilon},\quad a_\varepsilon:=m(E)\frac{H(\varepsilon^{-1})}{H(\varepsilon^{-\rho})}.
\end{equation*}
Recalling the probability law of big jumps \eqref{eq:law} we see that for
$\varepsilon$ small enough
\begin{equation}
\label{eq:k}
\mathbf{P}(F(0)\varepsilon J_k \notin \mathcal{G})=\mathbf{P}(\varepsilon J_k\in E)=\frac{\nu(E/\varepsilon)}{\beta_\varepsilon}
\end{equation}
whereas for any $\delta'>0$ 
\begin{equation}\label{eq:kk}
a_\varepsilon(1-\delta') \leq \frac{\nu(E/\varepsilon)}{\beta_\varepsilon}\leq 
a_\varepsilon(1+\delta').
\end{equation}
To obtain the final asymptotics we have to estimate carefully the perturbed the exit probabilities
$\mathbf{P}(F(y)\varepsilon J_1 \in \mathcal{G}^{-\varepsilon^\gamma})$ and $\mathbf{P}(F(y)\varepsilon J_1 \notin \mathcal{G})$ uniformly 
over $\|y\|\le \varepsilon^\gamma$. This is achieved with help of the continuity of the function $(y,z)\mapsto F(y)z$ both in $y$ and $z$.
Indeed, for any $\delta'>0$  we can choose $R>0$ big enough, such that the estimate 
\begin{equation*}
\mathbf{P}(\|\varepsilon J_k\|> R)\leq  \frac{\delta'}{4}   \frac{H(\varepsilon^{-1})}{H(\varepsilon^{-\rho})}
\end{equation*}
holds for $\varepsilon$ small.
Further,  the function $F(y)z$ is  uniformly continuous in $z$ in the ball $\|z\|\leq R$ and is continuous in
$y$ at the origin. Using the scaling property of the jump measure $\nu$ and the fact that the limiting measure $m$ has no atoms
 we show that
uniformly over $\|y\|\leq \varepsilon^\gamma$
\begin{equation*}
|\mathbf{P}(F(y)\varepsilon J_k \notin \mathcal{G}^{\pm\varepsilon^\gamma}, \|\varepsilon J_k\|\leq R) -
\mathbf{P}(F(0)\varepsilon J_k \notin \mathcal{G}, \|\varepsilon J_k\|\leq R )|  \leq \frac{\delta'}{4}\frac{H(\varepsilon^{-1})}{H(\varepsilon^{-\rho})},
\end{equation*}
and
\begin{equation*}
\mathbf{P}(F(0)\varepsilon J_k \notin \mathcal{G}) -
\mathbf{P}(F(0)\varepsilon J_k \notin \mathcal{G}, \|\varepsilon J_k\|\leq R ) \leq \frac{\delta'}{4}\frac{H(\varepsilon^{-1})}{H(\varepsilon^{-\rho})}.
\end{equation*}
Finally for any $\delta>0$ we choose $\delta'>0$ small enough to get the uniform estimates
\begin{equation*}
\begin{aligned}
& \inf_{\|y\|\le \varepsilon^{-\gamma}} \mathbf{P}(F(y)\varepsilon J_1\in \mathcal{G}^{-\varepsilon^\gamma} )
\geq 1-m(E)\frac{H(\varepsilon^{-1})}{H(\varepsilon^{-\rho})}(1+\delta)=1-a_\varepsilon (1+\delta) ,\\
&\inf_{\|y\|\le \varepsilon^{-\gamma}} \mathbf{P}(F(y)\varepsilon J_1\notin \mathcal{G} )
\geq m(E)\frac{H(\varepsilon^{-1})}{H(\varepsilon^{-\rho})}(1-\delta)=a_\varepsilon(1-\delta).
\end{aligned}
\end{equation*}
for $\varepsilon$ small enough.

Following the lines of the proof of \cite{ImkellerP-06a} and \cite{ImkPavSta-10}, for any $\delta>0$ and 
$\varepsilon$ small we can also obtain the multiplicative estimates for the Laplace transforms for any $u>-1$:
\begin{equation*}
\begin{aligned}
&\inf_{y\in \mathcal{G}^{-\varepsilon^\gamma}}\mathbf{E}\Big[e^{-u m(E) H(\varepsilon^{-1})T_1}\mathbb{I}\{X_{t}(y)\in \mathcal{G}^{-\varepsilon^\gamma}, t \in [0,T_{1}]\}\Big]
\geq \frac{1-a_\varepsilon (1+\delta)}{1+u a_\varepsilon},\\
&\inf_{y\in \mathcal{G}^{-\varepsilon^\gamma}}\mathbf{E}\Big[e^{-um(E) H(\varepsilon^{-1})T_1}\mathbb{I}\{X_{t}(y)\in \mathcal{G}^{-\varepsilon^\gamma}, t \in [0,T_{1})\} 
\mathbb{I}\{X_{T_1} \notin \mathcal{G}\}   \Big]
\geq 
 \frac{a_\varepsilon(1-\delta)}{1+u a_\varepsilon}.
\end{aligned}
\end{equation*}
Summing up the terms from \eqref{eq:tt} over  $k\geq 1$ yields the estimate
\begin{equation*}
\begin{aligned}
\mathbf{E} e^{-u m(E) H(\varepsilon^{-1}) \sigma_x(\varepsilon)}
&\geq \frac{a_\varepsilon(1-\delta)}{1+u a_\varepsilon}\sum_{k=1}^\infty 
\Big(\frac{1-a_\varepsilon(1+\delta)}{1 +u a_\varepsilon}\Big)^{k-1}\\
&=\frac{1-\delta}{1+u+\delta}\geq \frac{1}{1+u}-C(u)\delta
\end{aligned}
\end{equation*}
for some $C(u)>0$ and $\varepsilon$ small.

The upper bound for the Laplace transform is technically more involved since it additionally demands 
careful estimates of the probability to exit from the domain due to \textit{small jumps} during the inter-jump
intervals of the compound Poisson process $\eta$.
These estimates are obtained analogously to the one-dimensional and multi-dimensional cases studied in 
\cite{ImkellerP-06a,ImkPavSta-10} and finally
lead to the uniform convergence of the Laplace transform over $x\in \mathcal{G}^{-\varepsilon^\gamma}$ as $\varepsilon\to 0$.

\section{First exit time of the Stratonovich SDE}

Recalling the It\^o form of the Stratonovich SDE \eqref{eq:si}, we reduce the exit problem of $X^\circ$ to the It\^o case. 
Indeed, in the argument of the Section \ref{s:Ito} we have to take into account the Stratonovich correction term
$\frac{\varepsilon^2}{2}\int_0^t F'(X_s^\circ)F(X_s^\circ)d[Z,Z]^c$ which is 
a Lebesgue integral whose absolute value increases at most as $C\varepsilon^2t$ for some $C>0$.
It is clear that adding this term to the equation does not influence the estimates of the section \ref{ss:pert}. 
Thus the result follows immediately, and we obtain the same asymptotics of the first exit time as in 
the It\^o case.

\section{First exit time of the Marcus (canonical) SDE}

The analysis of the canonical Marcus SDE can also be reduced to the It\^o case.
As in Section \ref{s:Ito} above let us distinguish between big and small jumps of $Z$. Since the processes $\eta$ and 
$L=Z-\eta$ are independent, the Marcus equation can be rewritten in the It\^o form as 
\begin{equation*}
\begin{aligned}
X^\diamond_t=
x&-\int_0^t \nabla U(X^\diamond_s)\, ds +\varepsilon \int_0^t F(X^\diamond_{s-})\circ  dZ^c_s
+\varepsilon \int_0^t F(X^\diamond_{s-})\, d L^d_s\\
&+\sum_{s\le t} \Big(  \varphi^{\varepsilon \Delta L^d_s}(X^\diamond_{s-}) -X^\diamond_{s-}
+F(X^\diamond_{s-}) \varepsilon\Delta L^d_s\Big) 
+\sum_{s\le t} \Big(  \varphi^{\varepsilon \Delta \eta_s}(X^\diamond_{s-}) -X^\diamond_{s-}\Big).
\end{aligned}
\end{equation*}
Let us estimate the small jump correction term in the Marcus equation. 
The jumps of the process $\varepsilon L^d$ are bounded in absolute value by $\varepsilon^{1-\rho}$. The mapping
$u\mapsto \varphi_u^z$ is $C^2(\mathbb{R},\mathbb{R}^n)$, so the Talyor expansion yields
\begin{equation*}
\begin{aligned}
\varphi^z(x)=y(1,x;z)&=y(0,x;z)+\frac{d}{du}y(0,x;z)+\frac{1}{2}\frac{d^2}{du^2}y(\theta,x;z)\\
&=x+F(x)z+R(x,z), \quad \theta\in(0,1),
\end{aligned}
\end{equation*}
where
\begin{equation*}
|R^k(x,z)|\le    \frac{1}{2}\sup_{y\in\mathbb{R}^n}  \Big| \sum_{i,j=1}^mz_i z_j 
\sum_{l=1}^n  \frac{\partial F_{kj}(y)}{\partial x_l}F_{li}(y)\Big|,\quad  1\le k\le n.
\end{equation*}
 Since all $F_{li}$ are bounded with bounded derivatives, we obtain the estimate
\begin{equation*}
|\varphi^z(x)-x-F(x)z|\le C\|z\|^2
\end{equation*}
with some absolute constant $C>0$.
This leads to the inequality
\begin{equation}
\label{eq:z}
\Big|\sum_{s\le t} \Big(  \varphi^{\varepsilon \Delta L^d_s}(X^\diamond_{s-}) -X^\diamond_{s-}
+ F(X^\diamond_{s-})\varepsilon\Delta L^d_s    \Big)\Big|
\le C\sum_{s\le t}  \varepsilon^2\|\Delta L^d_s\|^2=
C\varepsilon^2 [L]^d_t,
\end{equation}
so that this summand is small due to Lemma \ref{l:q}. 
Thus we are again in the setting of the deterministic dynamical system $Y$ perturbed by a small noise process
\begin{equation}
\label{eq:n}
\varepsilon \int_0^t F(X^\diamond_{s-})\circ  dZ^c_s
+\varepsilon \int_0^t F(X^\diamond_{s-})\, d L^d_s
+\sum_{s\le t} \Big(  \varphi^{\varepsilon \Delta L^d_s}(X^\diamond_{s-}) -X^\diamond_{s-}
+F(X^\diamond_{s-}) \varepsilon\Delta L^d_s\Big)
\end{equation}
 and a big jump process $\varphi^{\varepsilon J_k}(X^\diamond_{\tau_k-})-X^\diamond_{\tau_k-}$ .

The arguments of the Section \ref{s:Ito} for the proof of the It\^o case can be applied to the Marcus 
canonical equation. First, due to the estimate \eqref{eq:z} and Lemmas \ref{l:q} and \ref{l:qq}
we obtain the exponential estimate of the Lemma \ref{l:p} for the solutions of the canonical Marcus SDE $X^\diamond$.

Then we again exploit the concept of the \textit{one big jump}, that is we show that
the exit from the domain $\mathcal{G}$ occurs with high probability at one of the 
jump times $\tau_k$. Just before the time $\tau_k$ the solution $X^\diamond$
stays in a small neighbourhood of the stable point, so the exit occurs if the jump $\varepsilon J_k$ is large enough, namely if 
$\varphi^{\varepsilon J_k}(X_{\tau_k-})\approx \varphi^{\varepsilon J_k}(0) \notin \mathcal{G}$. The events 
$\{\varepsilon J_1\notin E^\diamond\}=\{\varphi^{\varepsilon J_1}(0) \in \mathcal{G}\}$,\dots, 
$\{\varepsilon J_{k-1}\notin E^\diamond\}=\{\varphi^{\varepsilon J_{k-1}}(0)  \in \mathcal{G}\}$ and  
$\{\varepsilon J_k\in E^\diamond\}=\{\varphi^{\varepsilon J_k}(0) \notin \mathcal{G}\}$
 are independent and build up a geometric sequence of events. 
As in the It\^o case,  for any $\delta'>0$ their probabilities can be calculated in the limit of $\varepsilon\to 0$ as
\begin{equation*}
m(E^\diamond)\frac{H(\varepsilon^{-1})}{H(\varepsilon^{-\rho})}(1-\delta') \leq \mathbf{P}(\varphi^{\varepsilon J_k}(0) \notin \mathcal{G})\leq 
m(E^\diamond)\frac{H(\varepsilon^{-1})}{H(\varepsilon^{-\rho})}(1+\delta').
\end{equation*}
Fine estimates for the perturbed the exit probabilities
$\mathbf{P}(\varphi^{\varepsilon J_1}(y)\in \mathcal{G}^{-\varepsilon^\gamma})$ and $\mathbf{P}(\varphi^{\varepsilon J_1}(y)\notin \mathcal{G})$ are also 
obtained analogously to the It\^o case. Indeed, one can see the mapping $(y,z)\mapsto F(y)z$ as a particular 
case of the mapping  $(y,z)\mapsto \varphi^z(y)$ appearing in the Marcus equation. Again, $\varphi^z(x)$ is continuous in
$y$ at $y=0$ and is uniformly continuous w.r.t.\ $z$ in the ball $\|z\|\leq R$ with some $R$ big enough.  
Thus the argument of the Section \ref{s:Ito} can be repeated directly with $\varphi^z(x)$ instead of $F(y)z$ and $E^\diamond$
instead of $E$. Consequently for any 
$\delta>0$ we obtain the uniform estimates
\begin{equation*}
\begin{aligned}
& \inf_{\|y\|\le \varepsilon^{-\gamma}} \mathbf{P}(\varphi^{\varepsilon J_1}(y)\in \mathcal{G}^{-\varepsilon^\gamma}_1 )
\geq 1-m(E^\diamond)\frac{H(\varepsilon^{-1})}{H(\varepsilon^{-\rho})}(1+\delta),\\
&\inf_{\|y\|\le \varepsilon^{-\gamma}} \mathbf{P}(\varphi^{\varepsilon J_1}(y) \notin \mathcal{G} )
\geq m(E^\diamond)\frac{H(\varepsilon^{-1})}{H(\varepsilon^{-\rho})}(1-\delta)
\end{aligned}
\end{equation*}
for $\varepsilon$ small and hence also the estimates for the Laplace transform.

\section*{Acknowledgments}
The author is grateful to M.\ H\"ogele for interesting discussions and an 
anonymous referee for the careful reading of the manuscript and her/his valuable comments.


\end{document}